\documentclass[12pt,a4paper]{article}
\usepackage[utf8]{inputenc}
\usepackage[T1]{fontenc}
\usepackage{amsmath}
\usepackage{amsfonts}
\usepackage{amssymb}
\usepackage{graphicx}
\usepackage{amsthm}
\usepackage[all]{xy}
\usepackage{xcolor}
\usepackage{tikz}

\usepackage{amsthm}

\numberwithin{equation}{section}

\newtheorem{theorem}{Theorem}[section]

\newtheorem{definition}[theorem]{Definition}

\newtheorem{teorema}{Theorem}
\newtheorem{lema}{Lemma}

\author{Robson da Silva\footnote{supported by CNPq}, Almir Neto\footnote{supported by FAPEAM}, and Kelvin Souza\footnote{supported by FAPEAM}}
\title{Skew plane partitions according to the $m$th largest and $m$th smallest parts}
\date{}

\begin{document}
\maketitle

\begin{abstract}
We extend recent results by G. E. Andrews and G. Simay on the $m$th largest and $m$th smallest parts of a partition to the more general context of skew plane partitions. In order to do this, we introduce new objects called skew plane overpartitions.
\end{abstract}

\section{Introduction}

We remember that the \textit{Ferrers diagram} of a partition $\lambda = (\lambda_1, \ldots, \lambda_r)$ of $n$, $|\lambda| = n$, is an array of cells with $r$ left-justified rows and $\lambda_i$ cells in row $i$. We say that a partition $\mu = (\mu_1, \ldots, \mu_s)$ is \textit{contained} in $\lambda$, $\mu \subset \lambda$, if $s \leq r$ and $\mu_i \leq \lambda_i$ for $i=1, \ldots, s$. The \textit{skew shape $\lambda/\mu$}, where $\mu$ is contained in $\lambda$, consists of deleting those cells of the Ferrers diagram of $\lambda$ which are also cells of the Ferrers diagram of $\mu$.

Given a skew shape $\lambda/\mu$, a \textit{skew plane partition of shape $\lambda/\mu$} is a filling of $\lambda/\mu$ with nonnegative integers called parts such that rows and columns decrease weakly. These objects were introduced in \cite{Stanley} and studied by many authors, including \cite{Krattenthaler, Sagan}.

In order to define the new objects called skew plane overpartition, we need restrict ourselves to a certain set of skew plane partitions:

\begin{definition}
A skew plane partition of shape $\lambda/\mu$ is called square-free if there are not four repeated parts forming a square. 
\end{definition}

For example, consider the skew plane partition of shape $(5,5,3,2)/(2,2,1)$.  The first skew plane partition below is not square-free while the second one is square-free:
\begin{center}
\begin{tabular}{ccccc}
 &  & 3 & 3 & 1 \\
 &  & 3 & 3 & 1 \\
 & 4 & 2 & &\\
8 & 1 & & &
\end{tabular}
\hspace{1cm}
\begin{tabular}{ccccc}
 &	 & 7 & 3 & 1 \\
 &	 & 3 & 3 & 1 \\
 &	4 & 2 & &\\
8 & 	1 & & &
\end{tabular}
\end{center}

We remember from \cite{Corteel} that a \textit{plane overpatition} is a plane partition where (1) in each row the last occurrence of an integer can be overlined or not and all the other occurrences of this integer are not overlined and
(2) in each column the first occurrence of an integer can be overlined or not and all the other occurrence of this integer are overlined. Below we have one plane overpatition of 80:

\begin{center}
	\begin{tabular}{cccccc}
		7 & 6 & $\overline{6}$ & 5 & 3 & $\overline{3}$ \\
		$\overline{7}$ & 5 & 5 & 4 & &   \\
		6 & 4 & $\overline{4}$ & 3 & &   \\
		5 & $\overline{4}$ & 2 & 1 & &
	\end{tabular}
\end{center}

It is important to note that due to the two conditions above, every plane overpatition is a square-free plane partition if we remove the overlines. We are interested in plane overpatitions with one additional condition:
\begin{enumerate}
\item[(3)] if one integer is overlined then the last occurrence of this number in each row is overlined.
\end{enumerate}
As an example, in order to the above plane overpartition satisfying condition (3), each last appearance of the numbers 3, 4, 6, and 7 in each row must be overlined:
\begin{center}
\begin{tabular}{cccccc}
$\overline{7}$ & 6 & $\overline{6}$ & 5 & 3 & $\overline{3}$ \\
$\overline{7}$ & 5 & 5 & $\overline{4}$ & &   \\
$\overline{6}$ & 4 & $\overline{4}$ & $\overline{3}$ & &   \\
5 & $\overline{4}$ & 2 & 1 & &
\end{tabular}
\end{center} 

We define now the new objects we are going to work with.

\begin{definition}
Let $n$ and $m$ be integers such that $n \geq m \geq 0$, $\lambda$ and $\mu$ be partitions of $n$ and $m$, respectively, such that $\mu \subset\lambda$. A skew plane overpartition of shape $\lambda/\mu$ of $N$ is a square-free skew plane partition of shape $\lambda/\mu$ of $N$ satisfying conditions (1), (2), and (3). 
\end{definition}

According to this definition, every plane overpartition is a skew plane overpartition of shape $\lambda/\emptyset$, for some partition $\lambda$. Below we have two examples of skew plane overpartition of shape $(5,5,3,2)/(2,2,1)$:

\begin{center}
	\begin{tabular}{ccccc}
		&  & 8 & $\overline{3}$ & $\overline{1}$ \\
		&  & 3 & $\overline{3}$ & $\overline{1}$ \\
		& $\overline{3}$ & 2 & &\\
		$\overline{8}$ & $\overline{1}$ & & &
	\end{tabular}
	\hspace{1cm}
	\begin{tabular}{ccccc}
		&	 & 4 & $\overline{4}$ & $\overline{2}$ \\
		&	 & 2 & 2 & $\overline{2}$ \\
		&	$\overline{4}$ & $\overline{2}$ & &\\
		8 & $\overline{4}$ & & &
	\end{tabular}
\end{center}

The \textit{shadow} of a skew plane overpartition of shape $\lambda/\mu$ is the ordinary square-free skew plane partition of shape $\lambda/\mu$ obtained after removing the overlines. For example, the skew plane partition of shape $(5,4,2,2,1)/(3,2)$ of $28$
\begin{center}
	\begin{tabular}{ccccc}
		 &  &  & 2 & 1\\
		 & & 4 & 2 \\
		5 &4& &\\
		5&4 & &\\
		1
	\end{tabular}
\end{center} 
is the shadow of the following four skew plane overpartition of shape \linebreak $(5,4,2,2,1)/(3,2)$:

\begin{center}
	$ \begin{array}{ccccc}
	 &  & &\overline{2}& 1\\
	 & &\overline{4}&\overline{2}   \\
	\overline{5} &\overline{4}\\
	\overline{5}&\overline{4}\\
	1 	
	\end{array} \ \   
	\begin{array}{ccccc}
	 &  & &\overline{2}& \overline{1}\\
	& &\overline{4}&\overline{2}   \\
	\overline{5} &\overline{4}\\
	\overline{5}&\overline{4}\\
	\overline{1} \end{array}\ \  
	\begin{array}{ccccc}
	 &  & &\overline{2}& \overline{1}\\
	 & &\overline{4}&\overline{2}   \\
	\overline{5} &\overline{4}\\
	\overline{5}&\overline{4}\\
	1 \end{array}\ \  
	\begin{array}{ccccc}
	 & & &\overline{2}& {1}\\
	 & &\overline{4}&\overline{2}   \\
	\overline{5} &\overline{4}\\
	\overline{5}&\overline{4}\\
	\overline{1} \end{array}$
\end{center}

\

Let $PG_{j+l,k}(n)$ be the number of skew plane overpartitions $P$ of $n$ such that: $|\lambda| \leq n$, where $\lambda/\mu$ is the skew shape of $P$, $k$ is an overlined part, and exactly $j+l$ other parts (each larger than $k$) are overlined, where $l$ is the number of parts (each larger than $k$) repeated in the columns. 

Similarly, we denote by $PS_{j+l,k}(n)$ the number of skew plane overpartitions $P$ of $n$ such that: $|\lambda| \leq n$, where $\lambda/\mu$ is the skew shape of $P$, $k$ is an overlined part, and exactly $j+l$ other parts (each smaller than $k$) are overlined, where $l$ is the number of parts (each smaller than $k$) repeated in the columns.

It is important to note that those skew plane overpartitions enumerated by $PG_{j+l,k}(n)$ ($PS_{j+l,k}(n)$) do not have any repeated part smaller (larger) than $k$, since they would be overlined otherwise.

\

We wish to find formulas for computing:
\begin{itemize}
\item $pg_{m}(n,k)$, the number of square-free skew plane partitions $T$ of $n$ in which $k$ is the $m$th greatest part, i.e., there are exactly $m-1$ different integers larger than $k$, and $|\lambda| \leq n$, where $\lambda/\mu$ is the shape of $T$. For example, $pg_{2}(4,1)=27$. Indeed, the skew shapes $\lambda/\mu$ we shaw use satisfy $0 \leq |\mu| \leq |\lambda| \leq 4$. Then, we have 15 skew shapes that can be filled  either with 1 and 2 or with 1 and 3 in only one way, while there are 6 skew shapes, listed below, that can be filled  either with 1 and 2 or 1 and 3 in two different ways.
\begin{figure}[h!]
	\center
	\begin{tikzpicture}[scale=0.6][rounded corners, ultra thick]
	\draw (8.25,-1.0) node {\tiny $(3,1)/(1)$};
	\shade[top color=white,bottom color=white, draw=black] (8,0) rectangle +(0.5,0.5);
	\shade[top color=white,bottom color=white, draw=black] (8.5,0) rectangle +(0.5,0.5);
	\shade[top color=white,bottom color=white, draw=black] (7.5,-0.5) rectangle +(0.5,0.5);
	
	\draw (11.35,-1.0) node {\tiny $(3,1)/(2)$};
	\shade[top color=white,bottom color=white, draw=black] (11.5,0) rectangle +(0.5,0.5);
	\shade[top color=white,bottom color=white, draw=black] (10.5,-0.5) rectangle +(0.5,0.5);

	\draw (14.,-1.0) node {\tiny $(2,2)/(1)$};
	\shade[top color=white,bottom color=white, draw=black] (14,0) rectangle +(0.5,0.5);
	\shade[top color=white,bottom color=white, draw=black] (14,-0.5) rectangle +(0.5,0.5);
	\shade[top color=white,bottom color=white, draw=black] (13.5,-0.5) rectangle +(0.5,0.5);
	
	\draw (16.35,-1.5) node {\tiny $(2,1,1)/(1)$};
	\shade[top color=white,bottom color=white, draw=black] (16.5,0) rectangle +(0.5,0.5);
	\shade[top color=white,bottom color=white, draw=black] (16,-0.5) rectangle +(0.5,0.5);
	\shade[top color=white,bottom color=white, draw=black] (16,-1) rectangle +(0.5,0.5);
	
	\draw (19.15,-1.5) node {\tiny $(2,1,1)/(1,1)$};
	\shade[top color=white,bottom color=white, draw=black] (19,0) rectangle +(0.5,0.5);
	\shade[top color=white,bottom color=white, draw=black] (18.5,-1) rectangle +(0.5,0.5);
	
	\draw (21.5,-1.0) node {\tiny $(2,1)/(1)$};
	\shade[top color=white,bottom color=white, draw=black] (21.5,0) rectangle +(0.5,0.5);
	\shade[top color=white,bottom color=white, draw=black] (21,-0.5) rectangle +(0.5,0.5);
	\end{tikzpicture}
	\caption{The 6 skew shapes that can be filled in two different ways}
	\label{Fig2}
\end{figure}
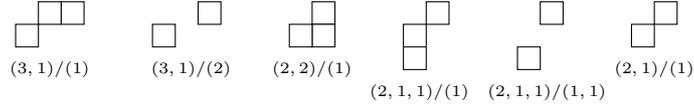

\item $ps_{m}(n,k)$, the number of square-free skew plane partitions $T$ of $n$ in which $k$ is the $m$th smallest part, i.e., there are exactly $m-1$ different integers smaller than $k$, and $|\lambda| \leq n$, where $\lambda/\mu$ is the shape of $T$. 
\end{itemize}

The theorems we shaw prove are:

\begin{teorema}
\textit{}
\begin{equation} pg_{m}(n,k)=  \sum\limits_{l \geq 0}\sum\limits_{j\geq 0}(-1)^{j+m-1-l} \binom{j}{m-1-l} PG_{j+l, k}(n).
\label{eq1}
\end{equation}
\label{th1}
\end{teorema}

\begin{teorema}
\textit{}
\begin{equation} ps_{m}(n,k)=\sum\limits_{l\geq 0}\sum\limits_{j\geq 0}(-1)^{j+m-1-l} \binom{j}{m-1-l} PS_{j+l, k}(n).
\label{eq2}
\end{equation}
\label{th2}
\end{teorema}

It is easy to see that Theorems 1 and 2 from \cite{Andrews1} are particular cases of the above theorems if we consider those overpartitions as skew plane overpartitions having only one row, which implies $l=0$.

\section{Proofs of Theorems 1 and 2}

It is easy to see that each square-free skew plane partition $P$ of $n$ having shape $\lambda/\mu$ and $D$ different parts is the shadow of $\binom{D-l}{j}$ skew plane overpartitions of $n$ having shape $\lambda/\mu$, in which exactly $j+l$ parts larger than $k$ are overlined, where $l$ is the number of parts that are repeated in the columns of the skew plane partition $P$ (each repeated part is counted only once). For example, in the square-free skew plane partition
$$\begin{array}{ccccc}
  & & 9 & 7 & 5\\
  & 4 & 3 & 2 & \\
5 & 3 & 3 & 1 & \\
5 & 3 & 2 & &
\end{array}$$
taking $k=2$, we have $D=5$ and $l=2$. Then, $0 \leq j \leq 3$, from where we see that there are $\binom{3}{0} = 1$, $\binom{3}{1} = 3$, $\binom{3}{2} = 3$, and $\binom{3}{3} = 1$ skew plane overpartitions having two, three, four, and five different overlined parts larger than $k$, respectively:
$$\begin{array}{ccccc}
& & 9 & 7 & \overline{5} \\
& 4 & \overline{3} & \overline{2} & \\
\overline{5} & 3 & \overline{3} & 1 & \\
\overline{5} & \overline{3} & \overline{2} & &
\end{array} \ \ \ \textnormal{} \ \ \ \ \ 
\begin{array}{ccccc}
& & \overline{9} & 7 & \overline{5} \\
& 4 & \overline{3} & \overline{2} & \\
\overline{5} & 3 & \overline{3} & {1} & \\
\overline{5} & \overline{3} & \overline{2} & &
\end{array} \ \ \ \textnormal{} \ \ \ \ \ 
\begin{array}{ccccc}
& & 9 & 7 & \overline{5} \\
& \overline{4} & \overline{3} & \overline{2} & \\
\overline{5} & 3 & \overline{3} & 1 & \\
\overline{5} & \overline{3} & \overline{2} & &
\end{array} \ \ \ \textnormal{} \ \ \ \ \ 
\begin{array}{ccccc}
& & 9 & \overline{7} & \overline{5} \\
& 4 & \overline{3} & \overline{2} & \\
\overline{5} & 3 & \overline{3} & 1 & \\
\overline{5} & \overline{3} & \overline{2} & &
\end{array}$$
$$\begin{array}{ccccc}
& & \overline{9} & 7 & \overline{5} \\
& \overline{4} & \overline{3} & \overline{2} & \\
\overline{5} & 3 & \overline{3} & {1} & \\
\overline{5} & \overline{3} & \overline{2} & &
\end{array} \ \ \ \textnormal{} \ \ \ \ \ 
\begin{array}{ccccc}
& & 9 & \overline{7} & \overline{5} \\
& \overline{4} & \overline{3} & \overline{2} & \\
\overline{5} & 3 & \overline{3} & {1} & \\
\overline{5} & \overline{3} & \overline{2} & &
\end{array} \ \ \ \textnormal{} \ \ \ \ \ 
\begin{array}{ccccc}
& & \overline{9} & \overline{7} & \overline{5} \\
& {4} & \overline{3} & \overline{2} & \\
\overline{5} & 3 & \overline{3} & 1 & \\
\overline{5} & \overline{3} & \overline{2} & &
\end{array} \ \ \ \textnormal{} \ \ \ \ \ 
\begin{array}{ccccc}
& & \overline{9} & \overline{7} & \overline{5} \\
& \overline{4} & \overline{3} & \overline{2} & \\
\overline{5} & 3 & \overline{3} & {1} & \\
\overline{5} & \overline{3} & \overline{2} & &
\end{array}$$

The following lemma is proved in \cite{Andrews1} and will be used in the proof of Theorems 1 and 2.
\begin{lema} Given nonnegative integers $R$ and $D$, we have
$$\sum\limits_{j\geq 0}^{}(-1)^{j+R} \dbinom{j}{R} \dbinom{D}{j} = \left\{\begin{array}{cc}
1, & \textnormal{if} \ D=R \\
0, & \textnormal{if} \ D\neq R. \\
\end{array}
\right.$$
\label{lemma1}
\end{lema}

We denote by $pg_{l,m}(n,k)$ ($ps_{l,m}(n,k)$) the number of square-free skew plane partitions enumerated by $pg_{m}(n,k)$ ($ps_{m}(n,k)$) having exactly $l$ repeated parts, each larger (smaller) than $k$, in the columns. 

\begin{proof}[Proof of Theorem \ref{th1}]
We proceed by examining the r.h.s of \eqref{eq1}. Let $P$ be an ordinary skew plane partition that is the shadow of some skew plane overpartitions enumerated by $PG_{j+l, k}(n)$. Thus $k$ must be a part of $P$ and there must be $D\geq j+l$ different parts of $P$ that are larger than $k$.
Now, we will find what is the contribution to the sum
\begin{equation}
\sum\limits_{j\geq 0}^{}(-1)^{j+m-1-l} \binom{j}{m-1-l} PG_{j+l, k}(n)
\label{eq3}
\end{equation}
of the skew plane overpartitions whose shadow is $P$.

From the observations in the beginning of this section, it follows that this contribution is $\sum_{j\geq 0}(-1)^{j+m-1-l}\binom{j}{m-1-l} \binom{D-l}{j}$. By Lemma \ref{lemma1} this contribution is either $0$, if $D \neq m-1$, or $1$, if $D=m-1$. In this latter case we have $j=m-1-l$ and therefore there exists exactly one skew plane overpartition counted by $PG_{j+l, k}(n)$, and dropping the overlines, we see that $P$ was a square-free skew plane partition having $k$ as the $m$th largest part. Thus in \eqref{eq3} the only contribution, of exactly $1$, come from skew plane overpartitions in one-to-one correspondence with square-free skew plane partitions in which $k$ is the $m$th largest part and having exactly $l$ repeated parts, each larger than $k$, in the columns.

Therefore,
$$pg_{m}(n,k)= \sum\limits_{l\geq 0}pg_{l,m}(n,k) = \sum\limits_{l\geq 0}\sum\limits_{j\geq 0}^{}(-1)^{j+m-1-l}\binom{j}{m-1-l} PG_{j+l, k}(n).$$
\end{proof}

The proof of Theorem \ref{th2} is similar to the above proof, except that now the other overlined parts, apart from $k$, are all smaller than $k$.


\section{Concluding Remarks}

In \cite{Andrews1} the authors provided an elegant recurrence for calculating the number of overpartitions with a given overlined greatest or smallest part and a given number of other overlined parts. They also present a generating function approach to their results. We hope to be able to get similar results for $PG_{j+l,k}(n)$ and $PS_{j+l,k}(n)$.

\

\noindent Robson da Silva  \\
Universidade Federal de S\~ao Paulo - Unifesp \\ 
S\~ao Jos\'e dos Campos - SP, 12247-014, Brazil \\
silva.robson@unifesp.br

\
 
\noindent Almir Neto \\
Universidade do Estado do Amazonas - UEA \\ 
Manaus - AM, 69050-010, Brazil \\
agneto@uea.edu.br

\

\noindent Kelvin Souza  \\
Universidade Federal do Amazonas - UFAM \\
Manaus - AM, 69077-000, Brazil \\
ksouzamath@yahoo.com.br


\begin{thebibliography}{99}
\bibitem{Andrews} G. Andrews, \textit{The  Theory  of  Partitions}. Cambridge Mathematical Library, Cambridge University Press, Cambridge, 1998.
\bibitem{Andrews1} G. Andrews and G. Simay, \textit{The $m$th Largest and $m$th Smallest Parts of a Partition}, \textit{Annal of Combinatorics}, accepted.
\bibitem{Corteel} S. Corteel, C. Savelief, and M. Vuleti{\'c}. \textit{Plane Overpartitions and Cylindric Partitions}. Journal of Combinatorial Theory Series A, 118(4), 1239 --1269, 2011.
\bibitem{Krattenthaler} C. Krattenthaler, \textit{Generating functions for plane partitions of a given shape}. Manuscripta Mathematica, 69, 173--201, 1990.
\bibitem{Sagan} B. E. Sagan, \textit{Combinatorial proofs of hook generating functions for skew plane partitions}. Theoretical Computer Science, 117, 273--287, 1993.
\bibitem{Stanley} R. P. Stanley, \textit{Theory and applications of plane partitions, Part 1}. Studies in Applied Mathematics, 50, 167--188, 269-279, 1971.
\end{thebibliography}
\end{document}